%% file: newdefs3.1.tex
\magnification=\magstep1

\def\to{\ \longrightarrow\ }

\def\nl{\hfill\break}

\input colordvi
\input epsf
\input cprfonts

\overfullrule=0pt

\font\npt=cmr9
\font\Bbb=msbm10

\font\secfont=cmbx10

\font\nam=cmr8
\font\aff=cmti8

\mathchardef\square="0\hexa03
\def\qed{\hfill$\square$\par\rm}
\def\np{\vfill\eject}
\def\boxing#1{\ \lower 3.5pt\vbox{\vskip 3.5pt\hrule \hbox{\strut\vrule
\ #1 \vrule} \hrule} }

\def\down#1{\ \lower 3.5pt\vbox{\vskip 3.5pt \hbox{\strut \ #1 \vrule} \hrule} }
\def\negdown#1{\ \lower 3.5pt\vbox{\vskip 3.5pt \hbox{\strut  \vrule \ #1 }\hrule} }

\def\item#1{\leavevmode\llap{#1\quad}}

\def\hneutral#1{\buildrel 0 \over #1}

\hsize=6.2 truein
\vsize=9 truein
\voffset=-0.5 in

\baselineskip=13 pt
\parskip=\baselineskip
 1

\parindent=0pt

\def\R{\hbox{\Bbb R}}

\def\oover#1{\vbox{\ialign{##\crcr
{\npt o}\crcr\noalign{\kern 1pt\nointerlineskip}
$\hfil\displaystyle{#1}\hfil$\crcr}}}

\newif \iftitlepage \titlepagetrue

\def\diagram#1{\global\advance\diagramnumber by 1
$$\epsfbox{newdefs3fig.\number\diagramnumber}$$
\centerline{\npt\bf Figure\ \the\secnum.\the\diagramnumber \npt\sl\ #1}\par}

\newcount\diagramnumber
\diagramnumber=0

\newcount\chapternumber \newcount\diagramnumber 
\newcount\questionnumber \newcount\chapternumber
\chapternumber=1 \diagramnumber=0 \questionnumber=0
\newcount\secnum \secnum=0\newcount\subsecnum \subsecnum=0
\newcount\defnum \defnum=1\newcount\theonum \theonum=0
\newcount\lemnum \lemnum=0

\newcount\subsecnum
\newcount\defnum
\def\section#1{
                \vskip 10 pt
                \advance\secnum by 1 
                \leftline{\secfont \the\secnum \rm\quad\bf #1}
                }

\def\subsection#1{
                \vskip 10 pt
                \advance\subsecnum by 1 
                 \leftline{\secfont \the\secnum.\the\subsecnum\ \rm\quad\bf \ #1}
                }

\def\definition#1{
                \advance\defnum by 1 
               \par \bf Definition
\the\secnum .\the\defnum  \it \ #1\rm \par
                }
   
    \def\theorem#1{
                \advance\theonum by 1 
                \par\bf Theorem  \the\secnum
.\the\theonum \sl\ #1\rm\ \par
               }

\def\lemma#1{
                \advance\lemnum by 1 
                \par\bf Lemma  \the\secnum
.\the\lemnum \sl\ #1\rm \ \par
                }

\def\corollary#1{
                \advance\defnum by 1 
                \par\bf Corollary  \the\secnum
.\the\defnum \sl\ #1\rm\ \par
               }

\def\exercise#1{
                \advance\defnum by 1 
                \par\bf Exercise  \the\secnum
.\the\defnum \sl\ #1  \par\rm
               }

\def\cite#1{
				\secfont [#1]
				\rm$\!\!\!$\nobreak\
}

\vglue 20 pt
\centerline{\secfont Alexander and Markov Theorems for Generalized Knots, I}
\medskip

\centerline{\nam Andrew Bartholomew}
\centerline{\nam Roger Fenn}
\centerline{\aff School of Mathematical Sciences, University of Sussex}
\centerline{\aff Falmer, Brighton, BN1 9RH, England}
\centerline{\aff e-mail: rogerf@sussex.ac.uk}
\vglue 20 pt


\centerline{\nam ABSTRACT}
\leftskip=0.5 in
\rightskip=0.5in

{\bf  In this paper we look at which Alexander and Markov theories can be defined for generalized knot theories. }

\leftskip=0 in
\rightskip=0 in
\baselineskip=13 pt
\parskip=\baselineskip

\parskip=\baselineskip
\smallskip

\section{INTRODUCTION}

The theorem that any knot can be represented by a braided diagram, i.e. the closure of a braid, was first proved by Alexander, \cite{A}, in 1923. The further result that any two braided diagrams representing the same knot can be joined by a sequence of Reidemeister moves in which the intervening diagrams are themselves braided was proved by Markov, \cite{Mar}, in 1935.

Since then there have been several re-proofs of these results, \cite{B}, \cite{Mor}, \cite{LR}, \cite{T}, \cite{V}, \cite{Y}. These all relate to classical knots. A paper on virtual and welded knots has been published by Kamada, \cite{Kam} and there is a paper on doodles by Gotin, \cite{G}.

In this paper, the second of a series on generalized knots, we shall prove analogues of the Alexander and Markov results for classes of generalized knots which we call regular and normal respectively. We do not know of any knot theories which are not regular. The examples of theories which are normal and so satisfy the hypothesis of the Markov type theorem include classical, virtual, welded, singular, virtual doodles and others. The theory of planar doodles is regular but not normal. So this theory satisfies an Alexander type theorem but we do not know if it also satisfies a Markov type theorem.

In a subsequent paper we will consider generalized braids, monoids, groups and the consequences of the results herein.

We would like to thank Colin Rourke for helpful discussions.

\section{Diagrams on the 2-sphere}
In this paper we will consider generalized knots and knot theories represented by diagrams on the oriented 2-sphere\footnote{$ ^1$}{There is no reason to restrict to spherical diagrams. In some sense a group is a knot with a presentation as a diagram.}  For more details, the interested reader should look at the first paper in the series \cite{F}. 

From now on we will usually drop the label `generalized' and not distinguish between knots and links. So generalized knots and links are called knots, generalized knot and link diagrams are called knot diagrams and so on.

A (spherical) knot diagram, usually denoted as, $K, L, M,\ldots$,  consists of the following.\nl
1. An immersion in general position of a compact closed 1-manifold into the sphere, $S^2=\R^2\cup\infty$. The image of one component of the 1-manifold is called a {\bf component} of the diagram.\nl
2. The double point crossings are labelled or tagged by a {\bf type} indicated by a roman letter, $a$ say.\footnote{$ ^2$}{In the first paper in this series the general crossings were indicated by roman letters $i$, $j$, $k$. Since this is a well used notation for the position of a braid crossing we have used $a$, $b$, $c$ instead.}
 The tags have a positive version, $a$, and a negative version, $\bar a$, which may not be different.\nl
3. Two diagrams, $K, L$, are considered the same if there is a homeomorphism of $S^2$ which takes one immersion to the other and preserves orientation and tags. We write, $K\cong L$ and call them {\bf isomorphic}.

\subsection{Seifert circles, graphs and trees}
In a diagram we can {\it smooth} crossings as follows, \cite{S}. Surround each crossing by an oblong neighbourhood called a {\bf crossing bridge}, labelled by the same tag as the crossing they replace, as in the figure 2.1.
\diagram{Smoothing a crossing}
The crossing bridge or bridge for short can now be drawn as a slightly thicker line. The diagram now becomes a {\bf Seifert graph}, consisting of a number of disjoint oriented {\bf Seifert circles}, which we will call {\bf cycles}, together with the bridges which join some pairs of cycles. The conected components of the cycles between the ends of the bridges are called the {\bf arcs} of the cycle.

Because a cycle is oriented, it is the boundary of a right hand disk and a left hand disk. 
A pair of cycles divide the sphere into two discs and a separating annulus. If both cycles are oriented in the same direction, i.e. are homologous cycles in the bounding annulus, then they are called a {\bf coherent} pair. Otherwise they are {\bf incoherent}.
Note that if two cycles are joined by a bridge then they are necessarily coherent. The annulus between a coherent pair is the intersection of a right and a left hand disk of the pair. A pair of cycles are {\bf adjacent} if there is a path from one to the other which is disjoint from the rest of the diagram.  

The components of the complement of the Seifert graph are called {\bf regions}. A region whose boundary is a single cycle is called a  {\bf polar region} and its boundary a {\bf polar cycle}. Note that the interior of polar regions are disjoint from the rest of the diagram if it is connected.

Let $h=h(K)$ be the number of incoherent pairs of the Seifert graph obtained from the diagram $K$. If $h=0$ the diagram is  {\bf braided}. 
 It is easy to see  that the number of polar regions is at least 2 and is only 2 if $K$ is  braided. If $K$ is braided then it is isomorphic to a diagram in which the cycles are circles of latitude and the bridges are arcs of longitude joining adjacent cycles.
 

 Let $K$ be a  knot diagram. If the diagram is connected we can define an oriented tree $T(K)$, see \cite{V}. The edges are in bijective correspondence with the cycles and two edges share a vertex if the associated cycles are adjacent. The flow from left disk to right disk across a cycle defines the orientation of the corresponding edge of the tree. 

The {\bf $n-$chain} tree, $C_n$,  is an interval divided into $n$ edges. In any tree, two edges $e_1, e_2$ can be connected by a unique $n-$chain in which  $e_1, e_2$ are the end edges. Let $s_1, s_2$ be cycles and let $e_1, e_2$ be the corresponding edges in $T(K)$. If the orientation of $e_2$ is the same as the orientation induced by $e_1$ along the chain then the cycles  $s_1, s_2$ are coherent and conversely. 
 
The following observation is at the heart of this paper.

\lemma{Suppose $K$ is a conected knot diagram on the 2-sphere. If $K$ is braided, then the cycles, which are all oriented coherently, are totally ordered by inclusion of their right(left) hand discs. In particular the tree, $T(K)$, is a coherently oriented chain. If $K$ is not braided, so $h(K)>0$, then there are a pair of adjacent incoherent cycles. } 

{\bf Proof:} If $K$ is braided the number of polar regions is  2 and $T(K)$ is a chain. Since $h(K)=0$, $T(K)$ is coherently oriented and provides a total order on the cycles of $K$. If $K$ is not braided the number of polar regions is greater than 2 and
$T(K)$ has at least 3 ends.  Let $e_1, e_2$ be the end edges of a chain in $T(K)$ which have opposite orientations. Then there will be edges $e'_1, e'_2$ in $T(K)$ which have opposite orientations and share a vertex. The corresponding cycles will be adjacent and incoherent.
\qed

\subsection{A worked example}
A knot diagram is shown below in figure 2.2 together with its Seifert graph and tree. The bridge types corresponding to the crossings are not indicated.
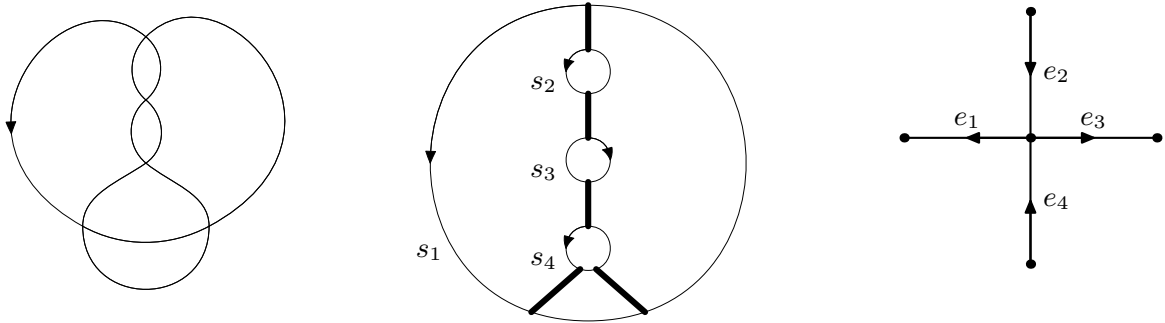
\diagram{Worked example: a knot diagram $K$, its Seifert graph and tree $T(K)$}
Let $s_1$ be the outermost cycle, $s_2$ the top one of the inner cycles,  $s_3$ the middle inner cycle and $s_4 $ the bottom one of the inner cycles. The edges $e_1,\ldots, e_4$ in the tree correspond to the cycles.

All 4 cycles are adjacent, $s_1$ is joined to $s_2$ and $s_4$ by bridges as is $s_3$. The cycle $s_1$ is  coherently oriented with $s_2$ and $s_4$ but not $s_3$. The cycles $s_2$ and $s_4$ are coherently oriented with $s_3$ but $s_2$ and $s_4$ are not coherently oriented.
So $h=2$.  All cycles are polar.

\section{The $R$ moves}

A {\bf Reidemeister} or $R$-move takes one diagram to another in one of the 4 ways indicated below in figure 3.3.

A generalized knot theory will define which of these moves is allowed and which are not.

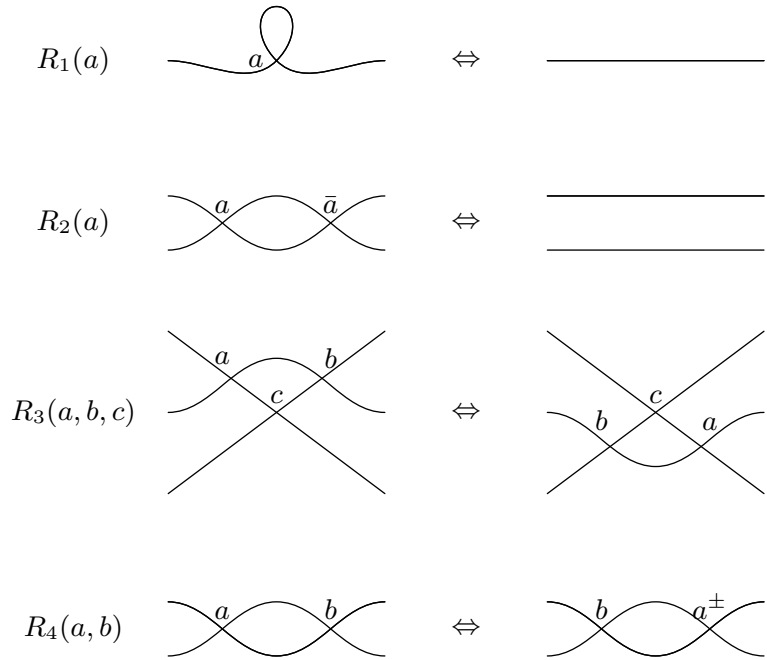
\diagram{The four $R$ moves}

\section{Orientation and the value of $h$}
In this section we look at possible orientations for the $R$ moves and how they impact on $h$.

An $R_1$ move creates/deletes a monogon. A creative move is denoted $R^+_1$ and a deleting move is denoted by $R^-_1$. The new monogon is a new cycle and is coherently oriented with its parent cycle. The monogon is called  a {\bf birth/death} disk. The value of $h$ is increased with the appearence of the monogon unless $h$ is originally 0 and the new cycle lies in one of the two polar regions. This special $R_1$ move  is called a {\bf Markov} move.

An $R_2$ move creates/deletes a bigon. As in $R_1$ moves they are noted $R^\pm_2$ according to creation or deletion. The bigon is a  birth/death disk as before.

An $R_2$ move can be divided into 2 kinds. If the arcs are oriented together, say from left to right then this is called an $R'_2$ move. These preserve the value of $h$.

If the arcs are oriented in opposite directions then this is called an $R_2''$ move. These change the value of $h$. If an $R_2''$ move involves the arcs from 2 distinct cycles then this is called a {\bf Vogel} or $V$ move.  These moves will play an important role in the subsequent proofs as we can see from the following lemma.

\lemma{A $V^+$ move which creates a bigon decreases $h$ by 1}
{\bf Proof: } The idea of the proof is to show that if $s, s'$ are two adjacent incoherent pairs of cycles then the $V$ move eliminates them and replaces them with a coherent pair. The details are as follows.

Assume that each cycle lies in the right disk of the other. Then the annulus between them is the intersection of their right hand disks. Let $l, a, r$ be the  number of cycles incoherent with $s$ which lie in the left hand disk, the annulus and the right hand disk but not 
$s'$ and not in the annulus, respectively. Define $l', a', r'$ similarly, but note that $a=a'$. Let $h(s)$ denote the number of cycles incoherently oriented with $s$. 
Then
$$h(s)=l+a+r+1 \hbox{ and } h(s')=l'+a+r'+1$$
Now do a $V$ move which eliminates $s, s'$ and introduces $c, c'$ where $c'$ is a polar cycle. Then
$$h(c)=l+l'+a \hbox{ and } h(c')=r+r'+a$$
Since $2h=\Sigma h(s)$, summed over all cycles, it follows that $h$ is reduced by 1.
\qed

If $R_2$ moves are always allowed then we call the theory {\bf regular}. Most importantly for regular theories, $V$ moves are allowed.

If the arcs of the $R_3$ move in figure 3.3 are all oriented from left to right, we specifically use the symbol 
$R'_3(a, b, c)$. If there is some crossing type $x$ such that $R'_3(x, \bar x, a)$  is true, then we say that $x$ {\bf dominates} $a$.  If there is some crossing type $x$ such that $R'_3(x, \bar x, a)$  and $R'_3(\bar x, x, a)$  are true for all $a$, then we say that $x$ {\bf dominates} the theory. A regular theory with a dominant $x$ is called {\bf normal}.

The other possibility for an $R_3$ move is that the central trigon is oriented. We denote this by $R''_3$.

Both versions of $R_3$ keep $h$ constant. 
An $R''_3$ move can be written in terms of $R'_3$ and $V$ moves, see \cite{F}.

The fourth move $R_4$ can be oriented in 2 ways. Either way does not alter $h$ and only changes tags.

\subsection{Finger and detour  moves}
Finger moves come in two sorts, $A$ and $B$ and are illustrated below in figure 4.4.  An $A$ finger move is applicable in a regular knot theory and consists of a series of $R_2$ moves in a line. A $B$ finger move is applicable in a normal theory with dominant tag $x$ and consists of $R_2$ moves  and an $R_3$ move using the dominant tag.

The tags, $y_1.\ldots, y_q$ are arbitary. The tags, $x_1, x_2, x_3, x_4$, take the values, $x, \bar x$ according to the orientations of the crossing arcs.

\diagram{$A$ and $B$ finger moves}

A subpath, $P$, of a diagram in a normal knot theory is said to be {\bf $x$ above} if the only crossings it meets are one or more of the two types illustrated below on the left of figure 4.5. The subpath, $P$, is drawn with a thicker line.

Similarly $P$  is said to be {\bf $x$ below} if the only crossings it meets are the ones on the right of figure 4.5.

\diagram{$x$ above and $x$ below paths}

The {\bf detour} move, see \cite{K}, is defined by the following lemma.

\lemma{{\bf (The detour move)} Let $P$ be an $x$ above/below subpath of a diagram $K$ in a normal theory and let $P'$ be a path with the same end points as $P$ which crosses $K$ in such a manner as to create an $x$ above/below path. Then the diagrams $K$ and $(K-P)\cup P'$ are related by a sequence of $R$ moves.}

{\bf Proof:} We can assume that $P\cup P'$ is the boundary of a bigon and argue by induction on the number of crossings inside. If there are none then an arc entering the bigon either leaves from the same side or the opposite side. If the former, then an innermost arc can be eliminated by an $R_2$ move. Eventually all the arcs cross from one side to the other. Then $P$ and $P'$ are isotopic.

If the bigon contains crossings then they can be eliminated one by one using $B$ type finger moves. \qed

\section{Examples of Knot Theories}
A {\bf knot theory} is defined by an $R$-graph. This is a graph with vertices, consisting of diagrams with an assortment of tags.  Two vertices (diagrams) are joined by an edge if there is an {\bf allowable} $R$-move between them. Note that if an $R$-move is not allowed then it is {\bf forbidden}. A knot from this theory is a component of the $R$-graph. \footnote{$ ^3$}{Note that the $R$-moves are invertible and so the $R$-graph generates a {\bf groupoid}.}

Some tags have specific properties due to their allowable $R$-moves either with themselves or with other tags. These properties are constant and define the knot theories with which they are involved. They have their own pictures with some overlap.

\subsection{Classical Knot Theory}

Classical or real crossings $r$ and $\bar r$ are indicated by a gap in the underarc and are illustrated as in figure 5.6.

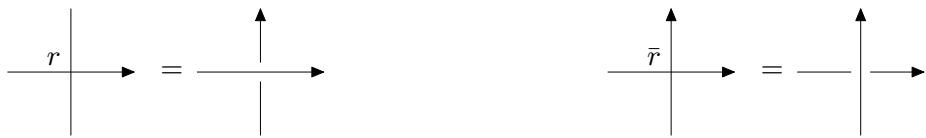
\diagram{Positive and Negative Classical Crossings}
They are either positive (right handed) or negative (left handed). Unless framings are considered $R_1$ and $R_2$ hold in a classical diagram. The allowed $R_3$ moves are defined by the fact that both $r$ and $\bar r$ dominate. This is a normal theory. Geometrically the allowed $R_3$ moves correspond to an arc of the diagram moving under or over a crossing respectively. So a detour move goes either over when the path is continuous or under when the path is broken at each crossing encountered.
\subsection{Virtual and welded knots}
\cite{K, FRR} The tags, $v$ and $w$ for virtual and welded crossings are both illustrated as a crossing with an enclosing circle as in figure 5.7.

\diagram{A vitual or welded crossing}

Both satisfy $R_1$ and $R_2$ and are involutive. So $v=\bar v$ and $w=\bar w$. 

The virtual crossing type $v$ dominates every other tag. So an arc of virtual crossings can move modulo its end points anywhere in the diagram by a detour move. Moreover no other tag can dominate $v$, so it is top of a sort of hierarchy.

The welded crossing is like a virtual crossing but also satisfies $R_3(\bar r, w, r)$, the {\bf first forbidden move} but not  $R_3( r, w,\bar r)$ the {\bf second forbidden move} so $r$ does not dominate.

\subsection{Free knots and doodles.}
The unadorned or flat crossing, $f$, can also stand for a free crossing, $F$, \cite{MI} or a doodle crossing, $d$, \cite{Kh, FT, BFKK}. All are involutive and satisfy $R_1$ and $R_2$. The crossings $f$ and $F$ dominate themselves but $d$ does not. This means that in a planar doodle diagram an arc cannot move past a crossing, momentarily creating a triple point.

A virtual doodle diagram has both $v$ and $d$ tags.\footnote{$ ^4$}{This is a generalization from the original definition given in \cite{FT}}

A free knot satisfies $R_4(F,v)$.

\subsection{Singular knots}

\cite{B2, FRR} Singular crossings are illustrated in figure 5.8.

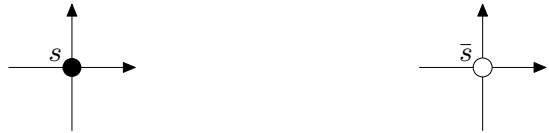
\diagram{Positive and Negative Singular Crossings}
The singular crossings $s$ and $\bar s$ do not satisfy $R_1$ but as usual $R_2$ is satisfied. They satisfy the fourth $R$-move, $R_4(r,s)$. In addition they do not dominate anything and are dominated by all other tags, so they are at the bottom of the heirarchy.

All of the above examples are normal. The original definition of doodles with only the $d$ tag is an example of a theory which is not normal.

\section{Alexander' theorem: an application of the Vogel move}

In this section we show that any diagram can be braided, (Alexander's theorem) provided Vogel moves are allowed.

\theorem{{\bf (The Braiding Process)} In a regular knot theory, a series of positive Vogel moves, one for each value of $h$,  can make all the Seifert cycles of a diagram coherently oriented. So $h$ becomes $0$, the Seifert cycles are nested and the diagram is braided.}
{\bf Proof:} Let $K$ denote the diagram and  suppose $h(K)>0$. By lemma 2.1 there is a pair of incoherently oriented cycles which are also adjacent. Push the arc of one along the path joining one to the other to make a $V^+$ move. This reduces $h$ by 1. If $h$ is still positive continue this process until $h=0$. \qed

\section{Generators for the $R$ moves}
In what follows, except in the proof of lemma 7.5, we will leave off the tags in diagrams and assume that any $R$ moves shown are allowed.

Recall that the $R$ moves which change $h$ are $R_1$, unless it is a Markov move, $R''_2$ moves involving arcs from the same cycle and $V$ moves.
In this section we will show that $R''_2$ moves which are not $V$ moves can be composed of $V$ moves and moves which do not change $h$.
If $h=0$ then the same is true for $R_1$ moves.

\lemma{Let $r:K\rightarrow L$ be a positive $R''_2$ move which is not a $V$ move in a regular theory. Then 
$r$ is a combination of $V$ moves, $R_1$ moves and moves which do not change $h$.}

{\bf Proof:} 
Figure 7.9 shows how the move $K\rightarrow L$ can be rewritten as a sequence of 2 $R_1$ moves, a positive $V$ move and a negative $V$ move. The up or down arrows indicate whether or not $h$ increases or decreases. \qed
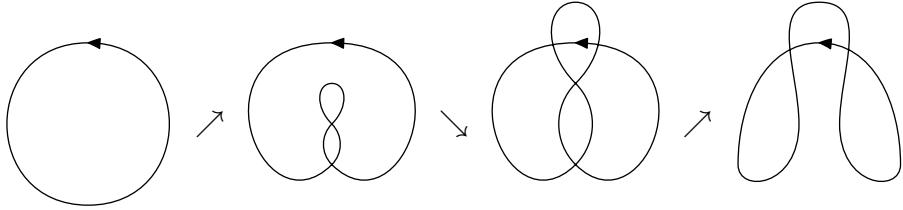
\diagram{Writing an $R''_2$ move in terms of other moves }

\lemma{Let $r:K\rightarrow L$ be a positive $R_1$ move in a normal theory. If $h(K)=0$ then $r$ can be written as a combination of $V$ moves and moves which do not change $h$.}

{\bf Proof:} 
Let the new crossing be tagged with $a$. Let $x$ be a dominant tag.

Since $h$ is raised, the birth disk must lie in an annulus defined by the parent cycle and one of the say $p>0$ cycles between the birth disk and a polar region. Then $K\rightarrow L$ can be written as a finger move sequence of $p$ $V$ moves tagged by $\bar x, x$ into the polar region, a Markov move tagged by $a$, $p$ $R'_3$ moves and a negative finger move as illustrated by figure 7.10. \qed
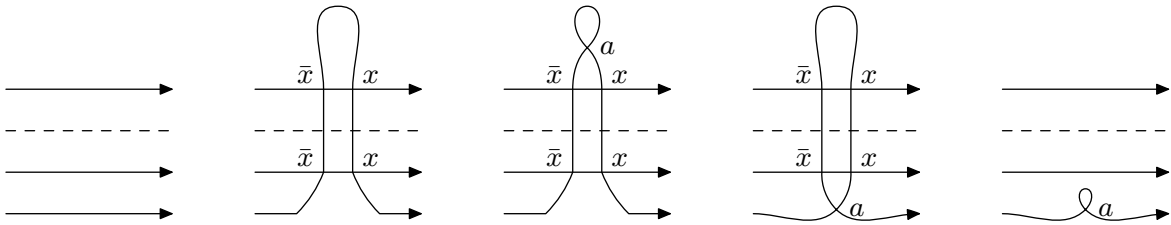
\diagram{Writing an $R_1$ move in terms of other moves}
This is the first time in the paper that normality has been needed.

\section{Markov's theorem}
In this final section, we will prove the result mentioned in the introduction. The idea behind this proof is again to use Vogel moves to lower the values of $h$ but their application will be more dynamic than in the proof of the generalized Alexander theorem since the elements in a sequence of $R$ moves will change during the application.

In view of the previous section we can assume that, in a normal theory, the only $R$ moves which change $h$ are the $V$ moves and 
non-Markov $R_1$ moves.

Firstly we will need some notation and definitions. A single $R$ move will be denoted by $K\rightarrow L$. A sequence of $R$ moves such as $K=K_1\rightarrow K_2\rightarrow \cdots \rightarrow K_{q-1}\rightarrow K_q=L$ of length $q$ will be denoted by $K\rightarrow\cdots\rightarrow L$. If the moves do not alter $h(K)$ we write 
$K\hneutral{\rightarrow}L$ or $K\rightarrow\hneutral{\cdots}\rightarrow L$ and refer to them as {\bf $h$-neutral} moves. 
A single $R$ move which raises (lowers) $h$ is denoted by $K\nearrow L$ ($K\searrow L$). The maximum value of $h$ in a sequence is denoted by $h(K\rightarrow\cdots\rightarrow L)=\max(h(K_i)), i=1,\ldots,q$. A sequence $K\nearrow L\searrow M$ is called a {\bf peak} and a sequence $K\searrow L\nearrow M$ is called a {\bf divot}.

We can now state the main theorem of this section.

\theorem{{\bf (Markov) }Suppose $K, L$ are braided and define the same knot in a normal theory, so they are related by a sequence $K\rightarrow\cdots\rightarrow L$. Then they are related by a sequence, 
$K\rightarrow\hneutral{\cdots}\rightarrow L$ in which all the intermediate diagrams are braided.}
The proof will follow from the following 2 lemmas.

\lemma{{\bf (The Peak Lemma)} Suppose $K\nearrow L\searrow M$ is a {peak} in a normal theory. Then one of the following is true.\nl
1. $K\cong M$,\nl
2. there is a sequence $K\rightarrow L_1\cdots\rightarrow L_q\nearrow M$ with $h(K\rightarrow L_1\cdots\rightarrow L_q\nearrow M) < h(L)$}

{\bf Proof:} Remember that we are only considering $V$ moves and non-Markov $R_1$ moves when changing $h$.

The proof is divided into a number of cases:\nl
i) both moves are $V$ moves\nl
ii) the increasing move is a $V$ move and the decreasing move is an $R_1$ move\nl
iii) the increasing move is an $R_1$ move and the decreasing move is an $R_1$ move\nl
iv) the increasing move is an $R_1$ move and the decreasing move is a $V$ move

Note that if $h(K)=0$ then by lemma 7.5 we do not need to consider cases iii) and iv).

The cases are further subdivided by the number of cycles involved, which can be 2, 3 or 4.

In case i), if the moves involve the same 2 arcs of a pair of cycles then 1. follows and $K\cong M$. If the tracks of the moves are disjoint then the moves can be interchanged and we have a divot which is outcome 2.

If there are three arcs of three cycles involved then the peak can be illustrated diagramatically by figure 8.11.

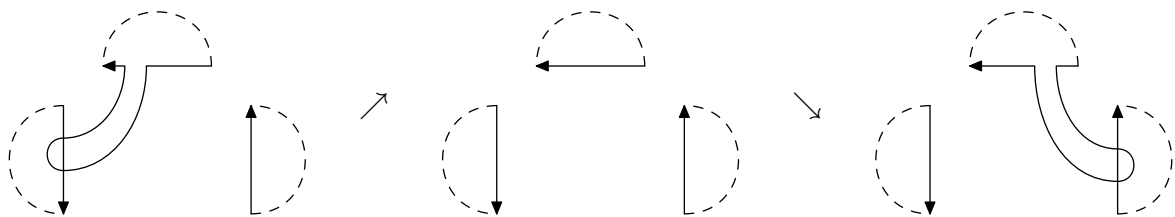
\diagram{A peak of $V$ moves with three arcs}

We now replace the middle $L$ by $L'$ as in figure 8.12.

\diagram{$L'$ the new $L$}

We now have a divot.

If there are four arcs of four cycles involved and the tracks of the $V$ moves cross, then this can be illustrated diagramatically by figure 8.13.

\diagram{The tracks of the $V$ moves cross}

We now replace the middle $L$ by $L'$ as in figure 8.14.

\diagram{$L'$ the new $L$}
To get to $L'$ from $K$ involves one $h$ neutral $R_2$ move and two positive $V$ moves. So $h$ is lowered by 2. To get from $L'$ to $M$ involves two negative $V$ moves and one $h$ neutral $R_2$ move. The peak is therefore replaced by outcome 2.

In case ii), the death disk of the $R_1$ move must be disjoint from the death disk of the increasing $V^{-}$ move and therefore the two moves can be interchanged, giving outcome 2.

In case iii), the birth and death disks are either equal or disjoint, which leads to outcome 1. or outcome 2. respectively.

In case iv), since we are only considering $h(K)>0$, by lemma 2.1 there is a positive $V$ move $r$ say whose track is disjoint from the birth disk and which lowers $h$ by 1. If we start with $r$, then do the peak and then $r^{-1}$ we will have outcome 2.\qed

\lemma{{\bf (The Transport Lemma)} Let $K\nearrow L\hneutral{\rightarrow} M$ be a sequence of two $R$ moves, the first of which raises and the second keeps  $h$ constant. Then one of the following is true,\nl
1. the moves can be interchanged, $K\hneutral{\rightarrow} L\nearrow M$\nl
2.  there is a sequence $K\rightarrow L_1\cdots\rightarrow L_q\nearrow M$ such that $h(K\rightarrow L_1\cdots\rightarrow L_q) < h(M)$}

{\bf Proof:} As with the Peak lemma there are several cases to consider.  The increasing move may be a $V^{-}$ move or, if $h(K)>0$, a non-Markov $R_1$ move.  For each of these moves, there are four possibilities for the $h$-neutral move:\nl
 i) a Markov $R_1$ move, \nl
 ii) an $R'_2$ move, \nl
 iii) an $R_3$ move, or \nl
 iv) an $R_4$ move.  
 
 Note that in case i) there is a regular neighbourhood of the $R_1$ move that is 
disjoint from the rest of the diagram and in cases ii), iii) and iv) there is a regular neighbourhood of the $h$-neutral move, within which
the move affects only the crossing bridges and not the cycles.

Assume initially that the first move is a $V$ move. Then, if its track does not interfere with the second move, they can be interchanged;  this will always be the situation in cases i), iii) and iv).  

So, assume that the other move is an $R'_2$ move involving 4 cycles and that their paths cross. This is a similar situation to the one described in figure 8.13 but the direction of the top path is reversed. 

We can repace $K\nearrow L\rightarrow M$ by the sequence $K\searrow L_1\nearrow L_2\rightarrow L_3\nearrow M$ as in the figure below.
\diagram{$K\searrow L_1\nearrow L_2\rightarrow L_3\nearrow M$}
The first two moves are $V$ moves, the next is the required $R'_2$ move and the last is another $V$ move. Since $h(K\searrow L_1\nearrow L_2\rightarrow L_3)<h(M)$ we are in outcome 2.

Suppose now that $K\nearrow L$ is a positive $R_1$ move. If the birth disk is disjoint from a regular neighbourhood of the $h$-neutral move, then we may interchange the order of the moves to obtain outcome 1.  If the birth disk intersects the regular neighbourhood, then it may be considered to lie wholly within it and, by lemma 2.1 there is a positive $V$ move $r$ say whose track is disjoint from the regular neighbourhood and which lowers $h$ by 1.  As in the proof of the peak lemma, we can therefore start with $r$, do the 2 moves and then do $r^{-1}$ and we will have outcome 2. 

This concludes the proof. \qed 

{\bf Proof of the generalized Markov:} Suppose the braided diagrams $K, L$ are related by a sequence $K\rightarrow\cdots\rightarrow L$. Consider a subsequence plateau $K'\nearrow K'' \rightarrow\hneutral{\cdots}\rightarrow L''\searrow L'$ where $h(K'' \rightarrow\hneutral{\cdots}\rightarrow L'')$ is maximal. By repeated applications of the transport lemma we can either eliminate the plateau or move the $h$ raising move to the right until a peak is formed. Then the peak lemma allows us to reduce the value of $h$. Eventually all the values of $h$ are reduced to zero and the theorem is proved. \qed

\section{Bibliography}
[A] Alexander, James (1923). "A lemma on a system of knotted curves". Proc. Natl. Acad. Sci. USA. 9: 93–95.

[BFKK] A. Bartholomew, R. Fenn, N. Kamada, S. Kamada "Doodles on Surfaces" to appear in JKTR

[B] J. Birman, "Braids, links and the mapping class groups", Annals of Math. Stud. 82, Princeton University Press, 1974.

[B2] J. Birman, "New points of view in knot theory",  Bull. Amer. Math. Soc. (N.S.) 28 (1993) 253-287

[F] Roger Fenn "Biquandles for Generalised Knot Theories" New Ideas in Low Dimensional Topology, pp. 79-103 (2015) Generalised  Series on Knots and Everything

[FRR] R Fenn, R Rimányi, C Rourke "The braid-permutation group" Topology 36 (1), 123-135

[FT] Roger Fenn and Paul Taylor," Introducing doodles", Topology of low-dimensional manifolds (Proc. Second Sussex Conf., Chelwood Gate, 1977), Lecture Notes in Math., vol. 722, Springer, Berlin, 1979, pp. 37–43.

[G] Konstantin Gotin  "Markov theorem for doodles on two-sphere" \nl arxiv.org/pdf/1807.05337

[Kam] Seiichi Kamada, "Braid presentation of virtual knots and welded knots",     Osaka J. Math. Volume 44, Number 2 (2007), 441-458.

[K] Kauffman, Louis H. (1999). "Virtual knot theory". European Journal of Combinatorics 20 (7): 663–690. 

[Kh] Mikhail Khovanov "Doodle groups" Trans. Amer. Math. Soc. 349 (1997), 2297-2315 

[LR]  Sofia Lambropoulou, Colin Rourke "Markov's theorem for 3-manifolds"
Topology and its Applications, 78 (1997) 95-112.

[Mar] A.A. Markov, "Über die freie Äquivalenz geschlossener Z\"opfe", Recueil Mathématique Moscou 1 (1935).

[Mor] H.R. Morton, "Threading knot diagrams", Math. Proc. Camb. Phil. Soc. 99 (1986), 247–260.

[MI]  Vassily Olegovich Manturov,  Denis Petrovich Ilyutko " Virtual Knots, the State of the Art"  World Scientific Series on Knots and Everything: Volume 51 (2012)

[N] Sam Nelson, "Unknotting virtual knots with Gauss diagram forbidden moves" J. Knot Theory Ramifications, 10, 931 (2001). 

[S] Seifert, H. "\"Uber das Geschlecht von Knoten" Mathematische Annalen (1935) pp. 571 - 592

[T] Paweł Traczyk, "A new proof of Markov's braid theorem" Banach Center Publications (1998) Volume: 42, Issue: 1, page 409-419  ISSN: 0137-6934

[V] P. Vogel, "Representation of links by braids: A new algorithm", Comment. Math. Helvetici 65 (1990), 104–113.

[Y] S. Yamada, "The minimal number of Seifert circles equals the braid index of a link" , Invent. Math. 89 (1987), 347–356.

\bye

\bye

\section{Generalised Braids}
Let $\cal K$ be a (generalized) knot theory with crossing types $T=\{ a,\ \bar a,\ b,\ \bar b\ldots\}$. Let $Y_n$ be the set of $n$  points, $Y_n=\{1, 2,\ldots{n}\}$ and let $X_m$ be a set of $m$  points in the interval $[0,1]$, $X_m=\{0=x_0<x_1<\cdots<x_{m-1}=1\}$.  

A braid diagram over $\cal K$ on $n$ strings consists of $n$ horizontal lines, $[0, 1]\times Y_n$ oriented from left to right and $m$ vertical bridge crossings, with labels from $T$, joining the $i$th line $[0, 1]\times \{i\}$ to the $i+1$ line $[0, 1]\times \{i+1\}$, where $i=i(j)$ and $x_j$ is the $x$-coordinate, $j=0, 1,\ldots, m-1$.

Braid diagrams multiply by concatonation in the usual way. There is an identity element with no bridges and if we agree that any sets $X_m$ are equivalent then the multiplication is associative.

Let $\sigma_i(a)$, be the braid diagram with one bridge labelled $a$  joining the $i$th line $[0, 1]\times \{i\}$ to the $i+1$ line $[0, 1]\times \{i+1\}$,  $i=1,2, \ldots, n-1$. These generate.
\diagram{ the generator $\sigma_i(a)$}
If we now add the relations, \nl
1. $\sigma_i(a)\sigma_j(b)=\sigma_j(b)\sigma_i(a)$ for $|i-j|>1$,\nl
2. $\sigma_i(a)\sigma_{i-1}(c)\sigma_i(b)=\sigma_{i-1}(b)\sigma_i(c)\sigma_{i-1}(a)$, $i=2, 3, \ldots, n-1$ whenever $(a,b)\succ c$, \nl
3. $\sigma_i(a)\sigma_i(b)=\sigma_i(b^\pm)\sigma_i(a)$ if $R_4(a,b)$.\nl
The braid diagrams form a monoid, $M_n({\cal K})$ generated by the $\sigma_i(a)$.

  The relations can be pictured in the following figures.
  \bye

%% file: cprfonts.tex
\def\hexnumber#1{\ifcase#1 0\or 1\or 2\or 3\or 4\or 5\or 6\or 7\or 8\or
 9\or A\or B\or C\or D\or E\or F\fi}
%
%
\font\twelvemsa=msam10 scaled 1200   
\font\tenmsa=msam10                  
\font\ninemsa=msam9            \font\sevenmsa=msam7
\font\sixmsa=msam6             \font\fivemsa=msam5
%
%
\newfam\msafam                 \textfont\msafam=\tenmsa
\scriptfont\msafam=\sevenmsa   \scriptscriptfont\msafam=\fivemsa
\edef\hexa{\hexnumber\msafam}        
\def\msa{\fam\msafam\tenmsa}         
%
%
\font\twelvemsb=msbm10 scaled 1200   
\font\tenmsb=msbm10                  
\font\ninemsb=msbm9            \font\sevenmsb=msbm7
\font\sixmsb=msbm6             \font\fivemsb=msbm5
%
\newfam\msbfam                 \textfont\msbfam=\tenmsb       
\scriptfont\msbfam=\sevenmsb   \scriptscriptfont\msbfam=\fivemsb
\edef\hexb{\hexnumber\msbfam}        
\def\msb{\fam\msbfam\tenmsb}         
%
%
\font\twelveeufm=eufm10 scaled 1200  
\font\teneufm=eufm10                 
\font\nineeufm=eufm9           \font\seveneufm=eufm7
\font\sixeufm=eufm6            \font\fiveeufm=eufm5
%
\newfam\eufmfam                \textfont\eufmfam=\teneufm
\scriptfont\eufmfam=\seveneufm \scriptscriptfont\eufmfam=\fiveeufm
\edef\hexf{\hexnumber\eufmfam}      
\def\frak{\fam\eufmfam\teneufm}     
%
%
%
\font\twelverm=cmr10 scaled 1200    
\font\ninerm=cmr9                   
\font\sixrm=cmr6   
%
\font\twelvei=cmmi10 scaled 1200    
\font\ninei=cmmi9                   
\font\sixi=cmmi6  
%
\font\twelvesy=cmsy10 scaled 1200   
\font\ninesy=cmsy9                  
\font\sixsy=cmsy6  
%
\font\twelvebf=cmbx10 scaled 1200   
\font\ninebf=cmbx9                  
\font\sixbf=cmbx6  
%
%
\font\twelveit=cmti10 scaled 1200   
\font\nineit=cmti9                  
%
\font\twelvesl=cmsl10 scaled 1200   
\font\ninesl=cmsl9                  
%
\font\twelvett=cmtt10 scaled 1200   
\font\ninett=cmtt9                  
%
%
%
%
\def\small{%
%
%
\textfont0=\ninerm \scriptfont0=\sixrm \scriptscriptfont0=\fiverm
\def\rm{\fam0\ninerm}        
%
%
\textfont1=\ninei \scriptfont1=\sixi \scriptscriptfont1=\fivei
%
%
\textfont2=\ninesy \scriptfont2=\sixsy \scriptscriptfont2=\fivesy
%
%
\textfont3=\tenex \scriptfont3=\tenex \scriptscriptfont3=\tenex
%
%
\textfont\bffam=\ninebf \scriptfont\bffam=\sixbf
\scriptscriptfont\bffam=\fivebf \def\bf{\fam\bffam\ninebf}%
%
%
\textfont\itfam=\nineit \def\it{\fam\itfam\nineit}%
\textfont\slfam=\ninesl \def\sl{\fam\slfam\ninesl}%
\textfont\ttfam=\ninett \def\tt{\fam\ttfam\ninett}%
%
%
%
\textfont\msafam=\ninemsa \scriptfont\msafam=\sixmsa
\scriptscriptfont\msafam=\fivemsa \def\msa{\fam\msafam\ninemsa}%
%
%
\textfont\msbfam=\ninemsb \scriptfont\msbfam=\sixmsb
\scriptscriptfont\msbfam=\fivemsb \def\msb{\fam\msbfam\ninemsb}%
%
%
\textfont\eufmfam=\nineeufm  \scriptfont\eufmfam=\sixeufm
\scriptscriptfont\eufmfam=\fiveeufm \def\frak{\fam\eufmfam\nineeufm}%
%
%
%
\normalbaselineskip=11pt
\setbox\strutbox=\hbox{\vrule height8pt depth3pt width0pt}%
%
%
\normalbaselines\rm}    
%
%
%
%
\def\large{%
\textfont0=\twelverm \scriptfont0=\ninerm \scriptscriptfont0=\sevenrm
\def\rm{\fam0\twelverm}%
\textfont1=\twelvei \scriptfont1=\ninei \scriptscriptfont1=\seveni
\textfont2=\twelvesy \scriptfont2=\ninesy \scriptscriptfont2=\sevensy
\textfont3=\tenex \scriptfont3=\tenex \scriptscriptfont3=\tenex
\textfont\bffam=\twelvebf \scriptfont\bffam=\ninebf
\scriptscriptfont\bffam=\sevenbf \def\bf{\fam\bffam\twelvebf}%
\textfont\itfam=\twelveit \def\it{\fam\itfam\twelveit}%
\textfont\slfam=\twelvesl \def\sl{\fam\slfam\twelvesl}%
\textfont\ttfam=\twelvett \def\tt{\fam\ttfam\twelvett}%
\textfont\msafam=\twelvemsa \scriptfont\msafam=\ninemsa
\scriptscriptfont\msafam=\sevenmsa \def\msa{\fam\msafam\twelvemsa}         
\textfont\msbfam=\twelvemsb \scriptfont\msbfam=\ninemsb
\scriptscriptfont\msbfam=\sevenmsb \def\msb{\fam\msbfam\twelvemsb}         
\textfont\eufmfam=\twelveeufm  \scriptfont\eufmfam=\nineeufm
\scriptscriptfont\eufmfam=\seveneufm \def\frak{\fam\eufmfam\teneufm}
\normalbaselineskip=15pt
\setbox\strutbox=\hbox{\vrule height11pt depth4pt width0pt}%
\normalbaselines\rm}%
%
\def\Bbb{\msb}

%

%
\mathchardef\plussquare="0\hexa01
\mathchardef\nge="3\hexb0B
\mathchardef\maltesecross="0\hexa7A
\mathchardef\del="0\hexf01
%

%